\documentstyle[12pt]{article}
\newtheorem{lem}{Lemma}[section]
\newtheorem{prop}{Proposition}

\newtheorem{defi}{Definition}[section]
\newtheorem{fact}{Fact}[section]
\newenvironment{proof}{{\bf Proof:}\newline}{\begin{flushright}$\Box$\end{flushright}}

\begin{document}
\title{One-Dimensional Maps and Poincar\'{e}\ Metric}
\author{Grzegorz \'{S}wia\c\negthinspace tek\\Institute for
Mathematical Sciences\\SUNY at Stony Brook\\Stony Brook, L.I., NY
11794\thanks{
This is a corrected and updated version of the Stony Brook 
preprint ``Bounded Distortion Properties of One-Dimensional Maps.''
Part of this work was done while the author was visiting
the Institute for Advanced Study in Princeton. 
Also, work supported in part by NSF grant 431-3604A.}}
\date{February 5, 1991}

\maketitle
\begin{abstract}
 
Invertible compositions of one-dimensional maps are studied which are
assumed
to include maps with non-positive Schwarzian derivative and others
whose sum of distortions is bounded. If the assumptions of the Koebe
principle 
hold,
we show that the joint distortion of the composition is bounded. On
the other hand, if all maps with possibly non-negative Schwarzian
derivative are almost linear-fractional and their nonlinearities tend
to cancel leaving only a small 
total, then they can all be replaced with affine maps with the same
domains and images and the resulting composition is a very good
approximation of the original one.  

These technical tools are then applied to prove a theorem about
critical circle  
maps. 
\end{abstract}

\section{Introduction}
\subsection{Review of results and techniques}
There are two ways of bounding the distortion of long compositions of
one-dimensional maps which appear so typically when we consider high
iterates of
a map. 

\subparagraph{Bounded nonlinearity.}
One is to use ``bounded nonlinearity" . The method goes back to
Denjoy.
 In the modern times, we think of nonlinearity of a function $f$ on a
one-
dimensional manifold as a form 
\[{\cal N}f=\frac{f''}{f'} dx\]
distributed
along the manifold and it turns out that the distortion of a high n-th
iterate on 
some interval is bounded by the integral of this form over the sum of
the images of this interval from the 0-th to (n-1)-st. There is a nice description of this
method with many applications to be found in \cite{dMvS}.

\subparagraph{``Koebe principle".}
If the map has critical points, its nonlinearity is non-integrable and
hence
no useful estimates can be obtained using the above mentioned method.
In this 
context, a new estimate was found in recent years. Instead of the
integrable
nonlinearity it uses negativity (positivity) of the Schwarzian
derivative.
It was first clearly stated  in \cite{guj}, though it seems that other
people 
had had similar ideas even before.  The Koebe principle gives pretty
good 
estimates, but the assumption of negative Schwarzian is unnervingly
strong to 
be made.

\paragraph{What we would like to know.}
There is another obvious observation, namely that any map has an
``integrable
nonlinearity" part and a ``negative Schwarzian" part. The Schwarzian
derivative 
must be negative in some neighborhood of each critical point, and
beyond the union
of these neighborhoods the nonlinearity is bounded. This observation
was made
and successfully used in a number of works. Estimates of the
distortion were
obtained, but they were typically estimates by large numbers.
Sometimes, it
is desirable to know also that the distortion is actually small. We
give this
kind of estimate in Section 2.

Another problem appears in conjunction with the study of universality,
notably
in the case of circle maps. It is widely believed that for circle
homeomorphisms whose rotation number is the golden mean their differentiable conjugacy class
depends only  on the type of their ``singularities"\footnote{meaning,
I guess,
points where the nonlinearity is infinite or undefined}. So far this
has been
proved in the situation of no singularities when it is the famous M.
Herman's
theorem (see \cite{Herman}, \cite{dMvS} and \cite{Skh}) and there is a
computer-assisted local argument in the case
of one cubic-type singularity \cite{Rand}. Herman's theorem implies,
and is
not too far from being equivalent to, that the first return maps on
small intervals tend
to be linear as the number of iterates involved grows. It means that
the distortions acquired by consecutive iterates of the map tend to cancel. We used the word cancel, because their sum with absolute values certainly {\em does not tend to zero.}

If we want to tackle the case when singularities do exist, we would at
least like  to know that we can asymptotically neglect the distortion coming from parts of the circle far from the singularities.  
The proof of this fact, called the ``pure singularity property''
occupies
the final sections of our work.

The emphasis of this paper is on technical problems, notably on the
methods
using the `` Poincar\'{e}\ model of the interval''. The possibility of
such an
approach was realized earlier and commented on by D. Sullivan 
(see \cite{miszczu}) and S. v. Strien. In the present paper, new
aspects
and applications of this technique are shown.

There are two main results of the paper: the Uniform Bounded
Distortion Lemma
in Section 2 and the Main Theorem in Section 4. The first result is a
tool
which I believe may be useful. The Main Theorem concerns universal
properties
of circle maps. The Uniform Bounded Distortion Lemma is not necessary
in order
to prove the Main Theorem, thus a reader who is only interested in the
Main 
Theorem has no need to advance beyond Lemma ~\ref{lem1} in Section 2.

This paper owes its inspiration to the graduate course taught 
by D. Sullivan in the fall of 1988. I also express my thanks to L.
Jonker,
A. Epstein, W. Pa\l uba and M. Samra whose keen remarks allowed me to
eliminate
a number of mistakes from the manuscript.

\subsection{Poincar\'{e}\  model of the interval}
We start with introducing the ``Poincar\'{e}\  model" of the interval.
If we are
given an interval $(\alpha,\delta)$ we can map its interior onto the
real line
by the map 
\[{\cal P}_{(\alpha, \delta)}(\gamma):=-log
(Cr(\alpha,\frac{\alpha+\delta}{2},\gamma,\delta))\; 
.\]

Here, $Cr(\alpha,\beta,\gamma,\delta)$ is the {\bf cross-ratio}
defined by
\[Cr(\alpha,\beta,\gamma,\delta)=\frac{(\beta-\alpha)(\delta-\gamma)}
{(\gamma-\alpha)(\delta-\beta)}\; .\]

For any interval $I$ let ${\cal A}_{I}$ be the affine map from $I$ to
$[0,1]$. It is
an easy observation that ${\cal P}_{I}={\cal P}_{[0,1]}\circ {\cal
A}_{I}$.
Another useful fact is that the distance between two points $x,y \in
(\alpha,
\delta)$  is equal to $|\log Cr(\alpha, x,y,\delta)|$. Thus, we may 
alternatively think of the Poincar\'{e}\  model as the interval
equipped with a new 
metric.

Whenever we have an map $\phi$ from an interval $I$ to an interval $J$
we
can consider $\overline{\phi}={\cal A}_{J}\circ\phi\circ {\cal
A}^{-1}_{I}$ mapping from
$[0,1]$ into itself. 
Then we define 
\[{\cal P}(\phi)={\cal P}_{[0,1]}\circ \overline{\phi}\circ {\cal
P}_{[0,1]}^{ -1}\;  .\] 

Thus we have defined the operator $\cal P$ which assigns to every map
from an interval to another interval its ``Poincar\'{e}\  model map''.
This definition depends on the choice of the domain and image of the
map. 

We will need to understand the action of this operator on the group of
orientation preserving self-homeomorphisms of $[0,1]$. The operator
then establishes an isomorphism with the group of
orientation-preserving
homeomorphisms of the real line. Linear-fractional maps of the
interval
become translations, maps with negative Schwarzian derivative are
mapped
to expandings of the line and maps with positive Schwarzian correspond
to contractions.

For every map on an interval $\phi$ we define its {\bf Poincar\'{e}\
distortion 
norm} (or simply distortion norm if there is no danger of confusion)
${\cal D}(\phi)$ by
\[{\cal D}(\phi)=\|{\cal P}(\phi)-id\|_{C^{0}}\; .\] 

Next, for any function $\phi$ defined on an interval $(x,y)$ we define

\begin{equation}
\rho(x,y;\phi) := \frac{|\phi((x,y))|}{|(x,y)|}\; .
\end{equation}
 
There is a correspondence between these quantities:
\begin{lem} \label {lem-1}
For any orientation preserving homeomorphism $\phi : (x,y)\rightarrow
J$ and 
$\gamma\in (x,y)$ with ${\cal P}_{(x,y)}(\gamma) = \gamma'$ 
\[ \log\frac{\rho(x,\gamma;\phi)}{\rho(\gamma,y;\phi} = {\cal
P}(\phi)(\gamma')
- \gamma' \; .\]
\end{lem}
\begin{proof}
Since we can pre- and postcompose $\phi$ with affine maps and that
will not 
change the quantities we are interested in, we may assume that $x=0,
y=1, 
J=(0,1) $. Then we simply compute 
\[ \log\frac{\rho(0,\gamma; \phi)}{\rho(\gamma,1;\phi)} = 
= \log\frac{(\phi(\gamma) - 0)(1-\gamma)}{(\gamma -
0)(1-\phi(\gamma))} = 
\log Cr(0,\phi(\gamma),\gamma,1) \; .\]

The absolute value of this quantity is the same as the Poincar\'{e}\
distance
between $\gamma$ and $\phi(\gamma)$ and the sign is correct provided
that 
$\phi$ preserves the orientation. The claim follows.
\end{proof}    

\paragraph{The meaning of the distortion norm $\cal D$.}
A more usual measure of distortion by an interval diffeomorphism $f$
whose
domain is $[\alpha, \delta]$ (i.e. it extends a little beyond 
$(\alpha,\delta)$ as a smooth map) is 
\begin{equation}\label{equ:-1,1}
 \sup\{ |\log\frac {f'(x)}{f'(y)}| : x,y \in [\alpha, \delta] \} \; .
\end{equation}

First, we notice that 
\[ \log\frac{f'(\alpha) |\delta-\alpha|}{|f(\alpha, \delta)|} = 
\lim_{\gamma \rightarrow\alpha} {\cal P}(f)(\gamma') - \gamma' \]
where $\gamma'$ is the image of $\gamma$ in the Poincar\'{e}\  model.
This follows
immediately from Lemma ~\ref{lem-1}. Since the analogous statement is
valid 
for $\delta$, we get that
\[ |\log\frac{f'(\alpha)}{f'(\delta)}| \leq 2 {\cal D}(f) \; .\]

Thus, the usual distortion norm given by formula ~\ref{equ:-1,1} is
bounded
by twice the supremum of the Poincar\'{e}\  distortion norms for all
restrictions 
of the map to subintervals.         

Fortunately, our future estimates of the Poincar\'{e}\  distortion
norm will have the
property that they are uniformly good for all restrictions to smaller 
intervals.

On the other hand, the norm given by formula ~\ref{equ:-1,1} is
obviously 
larger than the Poincar\'{e}\  distortion norm. 
\section{Uniform Bounded Distortion Lemma}

We will show an  estimate quite similar to the Koebe principle , save
that we 
will not require the
Schwarzian to be of a definite sign.  All we need is that the function
is a 
composition
of many functions, some of which have non-negative Schwarzian and the
joint 
distortion of others is bounded. This is what typically happens when
we 
consider a high iterate of a function.

\subsection{The formulation}
\paragraph{Standard compositions.}
We consider a function $f$ defined on an interval $(a,d)$ of the
following
form:
\begin{equation}\label{equ:2,1}
f:=\sigma_{m}\circ h_{m}\circ\cdots\circ\sigma_{1}\circ h_{1}\; ,
f_{0}=id\; .
\end{equation}

We will also  use the notation:
\[f_{k}:=\sigma_{k}\circ h_{k}\circ\cdots\circ\sigma_{1}\circ
h_{1},\;\; k\leq n\; .\]

All maps are defined on intervals and are order-preserving
homeomorphisms
onto the domain of the next map. Maps $\sigma_{i}$ are assumed to have
non-negative Schwarzian derivative.

Next, we define two distortion ``norms'' for any such composition.

\begin{itemize}
\item  The number $d_{1}$ is equal to\footnote{Caution! $d_{1}$ is
negative.}
\begin{equation}
\sum_{i=1}^{m}\inf\{0,\log\frac{\rho(\alpha,\beta;h_{i})\cdot\rho(\gamma,\delta;h_{i})}{\rho(\alpha,\delta;h_{i})\cdot\rho(\beta,\gamma;h_{i})}:\alpha<\beta<\gamma<\delta\in
f_{i}((a,d))\} 
\; .
\end{equation}
Since the argument of the logarithm represents the change of a
cross-ratio
under the map, the $i$-th term in the sum is zero provided that
$h_{i}$ has a 
non-negative Schwarzian derivative.
\item 
The number $d_{2}$ is
\begin{equation}
d_{2} :=  \sum_{i=1}^{m} {\cal D}(h_{i})
\end{equation}
\end{itemize}

The ``norm'' $d_{1}$ only gets closer to $0$ if we consider the
restriction
to a smaller interval. In order for $d_{2}$ to be uniform with respect
to
restrictions, it is sufficient to demand that the sum of 
``log-ratio of derivatives'' norms for maps $h_{i}$ be bounded.
 
In addition, we assume that
\[ \max\{|D(h_{i})| : 0 \leq i \leq n\} \leq \log 2 \; .\]

We will now see that the last requirement is only technical and can be
satisfied in each example.  More
precisely, if this condition initially is not satisfied, we can write
the same
function $f$ as another composition for which
the norms $d_{1}$ and $d_{2}$ stay the same, but the distortion norms
of 
functions $h_{i}$ become suitably small.

The way to do it is by rewritting all maps $h_{i}$ whose distortion is
too 
large as
compositions of many maps already with suitably small distortion
norms.
For each ${\cal P}(h_{i})$ we define the family
\[{\cal P}(h^{t}_{i}):=id+t({\cal P}(h_{i})-id)\; .\]
Obviously, $h^{1}_{i}=h_{i}$ and by dividing the interval $[0,1]$ into
sufficiently many subintervals we represent $h_{i}$ as a composition
of maps
with small distortion norms. The norms $d_{1}$ and $d_{2}$ will stay
the same.

Compositions satisfying these assumptions together with their norms
as defined above will be called {\bf standard compositions}.

\paragraph{Uniform Bounded Distortion Lemma.}
\subparagraph{The technical statement.}
\begin{em}

If a function $f$ is a standard composition of length $n$ defined on
an interval $(a,d)$, then for any interval $(b,c)\subset (a,d)$,
the distortion of $f$ on $(b,c)$ is bounded, namely:
\[{\cal D}(f_{|(b,c)})\leq\]
\[\leq Qd_{2}\exp(|d_{1}|)\min(1,\frac{c-b}{min(b-a,d-c)})+d_{2}+
2 |d_{1}+\log (Cr(a,b,c,d))|\]
where $Q$ is a constant quite independent of the composition.
\end{em}

\subparagraph{A simplified statement.}
\begin{em}
Suppose we have a standard composition defined on an interval $(a,d)$.
Then,
its Poincar\'{e}\  distortion on an interval $(b,c)$ is bounded by
\[d_{1} + d_{2} + K(d_{1},d_{2})|\log( Cr(a,b,c,d))|\]
where $K(d_{1},d_{2})$ is a constant depending only on $d_{1}$ and $d_{2}$ in 
a continuous fashion.
\end{em}

We leave it to the reader as an easy exercise to see that the
simplified 
version follows from the technical version.
\paragraph{A comment.}
We want to compare the classical Koebe priciple with our Uniform
Bounded 
Distortion Lemma and other estimates. One way to state the classical
Koebe
priciple is this:
\subparagraph{Koebe principle}
\begin{em}
If $g$ is a diffeomorphism defined on an interval $(a,d)$, and the
Schwarzian
derivative of $g$ is non-negative, then the nonlinearity coefficient
of 
$f$ is bounded, namely:
\[ |\frac{f''(x)}{f'(x)}| \leq \frac{2}{\min(x-a, d-x)} \; .\]
\end{em}

If we replace $\min(x-a, d-x)^{-1}$ with $1/(x-a) + 1/(d-x)$ and
integrate
from $y$ to $z$, we get
\[ |\log(f'(y)/f'(z))| < 2|\log Cr(a,y,z,d)| \; .\]

Thus, it becomes clear that the simplified version of the Uniform
Bounded
Distortion Lemma is a natural generalization of the Koebe principle.

As such, it is slightly stronger than estimates known so far (see
\cite{dMvS}
and \cite {blyu} for examples.) Those earlier estimates let us bound
the 
distortion by a uniform constant, while both the Koebe principle and
our 
lemma also give conditions for the distortion to be small (namely,
$c-b$ small
compared with the distance from $\{a,d\}$ and the ``$h$-contribution''
small.)
I know of one example, \cite{tanver}, when that makes a difference.
     
Finally, let us mention a strong recent result of \cite{miszczu} which
implies
that our $d_{1}$ ``norm'' can be controlled in terms of the Zygmund
norm of
$\bigcup h'_{i}$.\footnote{That is, in terms of the Zygmund norm of
the 
derivative in case when the standard composition is an iterate of the
map and
intermediate images of the domain are disjoint.}
   
\subsection{Proof of the Uniform Bounded Distortion Lemma.}
   The obvious approach to the proof of the Uniform Bounded Distortion
Lemma is by reducing the situation to the
Koebe principle. It will, however, require rearanging the order of
functions
$\sigma$ and $h$.  That would be an easy thing to do if both functions
had
the same domain:
\[\sigma\circ h = (\sigma\circ h\circ\sigma^{-1})\circ\sigma\] and
then we 
could regard the function in parentheses as a new $h$. Although the
interval 
maps usually do not have the same domains, their Poincar\'{e}\  model
maps
are all homeomorphisms of the whole line and so we play this trick
with them.

Then, two problems will appear:
first, whether after the rearangements the new functions $h$ will
preserve 
the bounded distortion properties; secondly whether we will indeed be
able to 
use the Koebe principle afterwards in order to bound the distortion of
Poincar\'{e}\  models of maps  $\sigma_{i}$.

\paragraph{The reshuffling procedure.}

The next lemma tells us what happens to the distortion norms of maps 
$h_{i}$ when we change the order of maps as described.
\begin{lem} \label{lem1}
Let $f$ be a standard composition. Then we can write ${\cal P}(f)$ as
the 
composition
\[\overline{h_{m}}\circ\ldots\circ\overline{h_{1}}\circ{\cal
P}(\sigma_{m})
\circ\ldots\circ{\cal P}(\sigma_{1}) \]
with 
\[\sum_{i=1}^{m}\sup\{|\overline{h_{i}}(\gamma)-\gamma|:\gamma\in
R\}\leq
\sum_{i=1}^{m}{\cal D}(h_{i}) \; .\]
\end{lem}
\begin{proof}
We start with rearanging the order of maps in the composition so as to
get
the contractions first.\footnote{A similar reshuffling, albeit in a
different 
context of complex quasiconformal extensions, was used by D. Sullivan
(see
\cite{miszczu}).}

We obtain the expression of the form:
\[\overline{h_{m}}\circ\overline{h_{m-1}}\circ\ldots\circ\overline{h_{1}}
\circ\tilde{\sigma}_{m}\circ\ldots\circ\tilde{\sigma}_{1}\; ,\]
where $\tilde{\sigma}_{i}$ means ${\cal P}(\sigma_{i})$ and
$\overline{h_{i}}$ is of the form \[\overline{\sigma}_{i}\circ {\cal
P}(h_{i})
\circ\overline{\sigma}_{i}^{\ -1}\] with 
$\overline{\sigma}_{i}$ being a certain composition of the functions 
$\tilde{\sigma}_{i}$ and therefore being a non-expanding map.

We compute: 
\[|\overline{h_{i}}(\gamma)-\gamma|=|\overline{\sigma}_{i}
\circ{\cal P}(h_{i})\circ\overline{\sigma}_{i}^{\ -1}(\gamma)-
\gamma|=|\overline{\sigma}_{i}
\circ{\cal P}(h_{i})\circ\overline{\sigma}_{i}^{\ -1}(\gamma) -
\overline{\sigma}_{i}\circ \overline{\sigma}_{i}^{\ -1}(\gamma)|\; .\]
Since $\overline{\sigma}_{i}$ is a contraction, the last expression
is not greater than
\[|{\cal P}(h_{i})\circ\overline{\sigma}_{i}^{\ -1}(\gamma)-
\overline{\sigma}_{i}^{\ -1}(\gamma)|\leq {\cal D}(h_{i})\; .\]
\end{proof}
   
We can apply Lemma~\ref{lem1} to the function $f$ restricted to the
interval
$(b,c)$ to find out that

\begin{equation}\label{1,1}
{\cal D}(f_{|(b,c)})\leq \sum_{i=0}^{m-1}{\cal
D}(h_{i+1|f_{i}((b,c))+})
\end{equation}
\[
\sup\{|{\cal P}(\sigma_{m|h_{m}\circ f_{m-1}(b,c)})\circ\ldots
\circ{\cal P}(\sigma_{1|h_{1}(b,c)})(\gamma)-\gamma|:\gamma\in R\}\; .
\]

The first sum in the inequality ~\ref{1,1} can be sufficiently sharply
bounded by $d_{2}$.  
What remains is to bound the distortion of maps with non-negative
Schwarzian
in the second term.
Naturally, we are going to achieve that using the Koebe principle,
however
we must be careful about the domains.

\paragraph{The issue of domains.}
We proceed as follows:
A map $g_{i}$ is defined to be the linear fractional map that agrees
with the map $h_{i}$ on the points $f_{i-1}(a),f_{i-1}(b),f_{i-1}(c)$.
At the same time, we assume without loss of generality that $b-a\leq
d-c$.
We further consider maps $F$ and $F_{i}$ defined in the analogous way
to
$f$ and $f_{i}$, save that the functions $h_{i}$ are replaced by
$g_{i}$. The 
function $F$ has non-negative Schwarzian, but in order to use the
Koebe 
principle we need to prove that it is defined on sufficiently big
neighborhood
of $(b,c)$. It is certainly defined on $(a,c)$ and we want to choose a
point $d'>c$ so that 
$F$ is defined on $(a,d')$ as well. This problem is solved by the next
lemma:
\begin{lem}\label{lem-2}
If we choose $d'$ so that
\[\frac{(d'-c)(b-a)}{(c-b)(d'-a)}= \exp(d_{1})\cdot\frac{(d-c)(b-a)}
{(c-b)(d-a)}\; ,\]
then $F$ is defined on $(a,d')$.      
\end{lem}
\begin{proof} 
 To check that we need to ensure that 
\[g_{i}(F_{i-1}(d'))<h_{i}(f_{i-1}(d))\] holds 
for every $i\leq m$. To prove that we will show that  
\begin{equation}\label{7,1}
\frac{(g_{i}(F_{i-1}(d'))-g_{i}(F_{i-1}(c)))\cdot 
(g_{i}(F_{i-1}(b))-g_{i}(F_{i-1}(a)))}{(g_{i}(F_{i-1}(c))-g_{i}(F_{i-1}(b)))
\cdot (g_{i}(F_{i-1}(d'))-g_{i}(F_{i-1}(a)))}\] \[\leq 
\frac{(h_{i}(f_{i-1}(d))-
h_{i}(f_{i-1}(c)))\cdot (h_{i}(f_{i-1}(b))-h_{i}(f_{i-1}(a)))}{(
h_{i}(f_{i-1}(c))-h_{i}(f_{i-1}(b)))\cdot
(h_{i}(f_{i-1}(d))-h_{i}(f_{i-1}(a)))} \; .
\end{equation}

The two complicated ratios above represent values of a certain
cross-ratio on
the images of points $a,b,c,d'$ and $a,b,c,d$ respectively.
Whenever we apply an order preserving homeomorphism $\phi$ to points
$\alpha<\beta<\gamma<\delta'<\delta$ the changes imposed on these
cross-ratios
will be related in the following way:
\[ \frac{\rho(\gamma,\delta';\phi)\cdot\rho(\alpha,\beta;\phi)}{
\rho(\beta,\gamma;\phi)\cdot\rho(\alpha,\delta';\phi)} :
\frac{\rho(\gamma,
\delta;\phi)\cdot
\rho(\alpha,\beta;\phi)}{\rho(\beta,\gamma;\phi)\cdot
\rho(\alpha,\delta;\phi)} =
\frac{\rho(\gamma,\delta';\phi)\cdot\rho(\alpha,
\delta;\phi)}{\rho(\alpha,\delta';\phi)\cdot\rho(\gamma,\delta;\phi)}\;
.\]
   The last expression is distortion of some kind of cross-ratio.
Using the
method from \cite{preston} we can show that if $\phi$ has a
non-negative 
Schwarzian, this number is not greater than $1$. That means that
provided
(\ref{7,1}) holds with some $i$, the subsequent application of
$\sigma_{i}$
will not increase the ratio of the left-hand side of (\ref{7,1}) to
the 
right-hand side. Next, we will apply $g_{i+1}$ to the points
$F_{i}(a),F_{i}(b),F_{i}(c),F_{i}(d')$ which is not going to change the cross-ratio at all; and
$h_{i+1}$ that will be applied to the points
$f_{i}(a),f_{i}(b),f_{i}(c),f_{i}(d)$ can decrease the cross-ratio by some factor $\zeta$. As this reasoning
shows, the ratio between the left-hand side of (\ref{7,1}) to the
right-hand side, will grow by no more than $\zeta$ as we pass from $i$ to $i+1$. Our assumptions imply that the total growth as we pass from $i=1$ to $i=n-1$ will be no 
more than $\exp(-d_{1})$.

So, the lemma follows.
\end{proof}

From now on, we take $(a,d')$ as specified in Lemma ~\ref{lem-2} to be
the
default domain of $F$.
 
By the Koebe principle we get the estimate for the nonlinearity
coefficient
of $F$, i.e. $nF := f''/f'$.

\[nF\leq
\frac{2}{min\{\gamma-a,d'-\gamma\}}<\frac{2}{\gamma-a}+\frac{2}{d'-\gamma}\; .\]
Integrating this inequality over $(b,c)$ we get
\[\sup\{|\log(F'(\gamma))-\log(F'(\beta))|:\gamma,\beta\in (b,c)\}<\]
\[-2\log(Cr(a,b,c,d'))=-2d_{1}-2\log(Cr(a,b,c,d) \; .\]
The quantity on the left-hand side of the last inequality is certainly
not 
less than ${\cal D}(F_{|(b,c)})$.

Thus,

\begin{equation}\label{2,1}
\sup\{|{\cal P}(\sigma_{m|(b,c)})\circ\ldots\circ{\cal
P}(\sigma_{1|(b,c)})
(\gamma)-\gamma|:\gamma\in R\}\leq
\end{equation}
\[ -2d_{1}-2\log(Cr(a,b,c,d)) + \sum_{i=1}^{m}{\cal
D}(g_{i|F_{i-1}(b,c)})\; 
.\]
 
\paragraph{The distortion of the linear-fractional part.}
The first two terms on the right-hand side of the inequality
~\ref{2,1}
can also be found the statement of the Uniform Bounded Distortion
Lemma.
What  remains to be calculated is the last sum. 

\begin{lem}\label{lem2}
Let $g$ be a linear fractional map defined on an interval $(a,c)$,
points $b$
and $\lambda$ belong to $(a,c)$ and satisfy $b<\lambda$. Then
\[\log\rho(b,\lambda;g)-\log\rho(\lambda,c;g)=\log\rho(a,b;g)-\log\rho(a,c;g)
\; .\] 
\end{lem}
\begin{proof}
This is a straightforward computation:
\[\log\rho(b,\lambda;g)-\log\rho(\lambda,c;g)=\log\frac{\rho(b,\lambda;g)}
{\rho(a,\lambda;g)}-\log\frac{\rho(\lambda,c;g)}{\rho(a,\lambda;g)}\]
\[=\log\frac{\rho(b,\lambda;g)\cdot\rho(a,c;g)}{\rho(a,\lambda;g)\cdot\rho(b,c;g)}+\log\frac{\rho(b,c;g)}{\rho(a,c;g)}+({\cal
P}(g)(P(\lambda))-P(\lambda))
\; .\]
We now notice two facts: that the first term is $0$, because it is
logarithm of
the distortion of some cross-ratio and $g$ is linear-fractional; next,
that we
can replace $P(\lambda)$ with any other point, provided ${\cal P}(g)$
is an
isometry. Therefore, the last expression is equal to
\[\log\frac{\rho(b,c;g)}{\rho(a,c;g)}+({\cal P}(g)(P(b))-P(b))=
\log\frac{\rho(b,c;g)}{\rho(a,c;g)}+\log\frac{\rho(a,b;g)}{\rho(b,c;g}=\]
\[\log\rho(a,b;g)-\log\rho(a,c;g)\]
\end{proof}

\paragraph{Final estimates.}
We are now ready to conclude the proof. From equation \ref{2,1} and 
Lemma~\ref{lem2} we get
\begin{equation}\label{4,1}
\sup\{|{\cal P}(\sigma_{m|h_{m}\circ
f_{m-1}(b,c)})\circ\ldots\circ{\cal P}(\sigma_{1|h_{1}(b,c)})
(\gamma)-\gamma|:\gamma\in R\}
\end{equation}
\[\leq -2\log(Cr(a,b,c,d'))-2d_{1}+\sum_{i=1}^{m}
|\log\frac{\rho(F_{i-1}(a),F_{i-1}(b);g_{i})}{\rho(F_{i-1}(a),F_{i-1}(c);
g_{i})}|\; .\]

We will prove a simple computational lemma that will enable us to
evaluate
the last sum.
\begin{lem}\label{lem3}
Suppose that a function $\phi$ is defined on an interval $(a,c)$ with 
${\cal D}(\phi)\leq \log 2$ and
a point $b\in (a,c)$. Then 
\[|\log\frac{\rho(a,b;\phi)}{\rho(a,c;\phi)}|\leq Q\cdot {\cal
D}(\phi)
\min(1,\frac{c-b}{b-a}) \; .\]
where $Q$ is a uniform constant. 
\end{lem}
\begin{proof}
\[|\log\frac{\rho(a,c;\phi)}{\rho(a,b;\phi)}|=|\log(\frac{\phi(c)-\phi(a)}
{\phi(b)-\phi(a)}:\frac{c-a}{b-a})|=|\log\frac{1+\frac{\phi(c)-\phi(b)}
{\phi(b)-\phi(a)}}{1+\frac{c-b}{b-a}}|\]
\[=|\log\frac{1+\exp(\log\rho(b,c;\phi)-\log\rho(b,a;\phi))\cdot
\frac{c-b}{b-a}}
{1+\frac{c-b}{b-a}}|\]
\begin{equation}\label{equ:7,1}
\leq \log \frac{1+\exp({\cal
D}(\phi))\cdot\frac{c-b}{b-a}}{1+\frac{c-b}{b-a}}
\end{equation}

Our final estimate will depend on how we bound the argument of the
logarithm
in Formula \ref{equ:7,1}.
One possible estimate is
\[\frac{1+\exp({\cal
D}(\phi))\cdot\frac{c-b}{b-a}}{1+\frac{c-b}{b-a}}\leq
\exp({\cal D}(\phi)) < 2 {\cal D}(\phi)\; .\]

If $\frac{c-b}{b-a}$ is small, we can expand formula \ref{equ:7,1}
into a 
series, and get an estimate proportional to $\frac{c-b}{b-a}$.
\end{proof}

We will apply Lemma~\ref{lem3} in the situation when $\phi:=g_{i}$,
$a:=F_{i-1}(a),b:=F_{i-1}(b), and c:=F_{i-1}(c)$.

We first note that ${\cal D}(g_{i})\leq {\cal D}(h_{i})$. Then we
recall
that the distortion norms of maps $h_{i}$ can be assumed to be as
small as
we want, in particular could be less than $\log 2$. Thus the
assumptions 
of Lemma ~\ref{lem3} are satisfied. 

It enables us to bound the second sum in Formula ~\ref{4,1} by

\begin{equation}\label{5,1}
Q\cdot \sum_{i=1}^{m} {\cal
D}(h_{i})\min(1,\frac{f_{i-1}(c)-f_{i-1}(b)}
{f_{i-1}(b)-f_{i-1}(a)})
\end{equation}
What we need is to estimate the ratios in (\ref{5,1}) in  terms of
$\frac{c-b}
{b-a}$.To do that, we consider a cross-ratio
$CR(\alpha,\beta,\gamma,\delta)$
defined by 
\[CR(\alpha,\beta,\gamma,\delta)=\frac{(\delta-\gamma)\cdot(\beta-\alpha)}
{(\gamma-\beta)\cdot(\delta-\alpha)}\; .\]

\[\frac{f_{i}(b)-f_{i}(a)}{f_{i}(c)-f_{i}(b)}>CR(f_{i}(a),f_{i}(b),f_{i}(c),
f_{i}(d))>\exp(-d_{1})\cdot CR(a,b,c,d)\]
\[>\exp(-d_{1})\frac{b-a}{c-b}\cdot \frac{d-c}{b-a+(d-c)+(c-b)}\geq
\exp(-d_{1})\frac{b-a}{c-b}\cdot\frac{1}{2+\frac{c-b}{b-a}} \]
hence
\[\frac{f_{i}(c)-f_{i}(b)}{f_{i}(b)-f_{i}(a)}<2\exp(-d_{1})\cdot
\frac{c-b}{b-a}(1+\frac{c-b}{b-a})\]
Thus, if $\frac{c-b}{b-a}<1$, we get
\[\frac{f_{i}(c)-f_{i}(b)}{f_{i}(b)-f_{i}(a)}<4\exp(-d_{1})\frac{c-b}{b-a}\;
.\]
We can finally bound the expression (\ref{5,1}) by
\[4Qd_{2}\exp(-d_{1})\min(1,\frac{c-b}{b-a})\]
which allows us to estimate the whole (\ref{4,1}) by
\begin{equation}\label{6,1}
-2\log(Cr(a,b,c,d)-2d_{1}+4Qd_{2}\exp(-d_{1})\min(1,\frac{c-b}{b-a})\;
.           %
\end{equation}

This concludes the proof of the Uniform Bounded Distortion Lemma.

\section{Functions $h_{i}$ with cancelling distortions}
In our formulation of the Uniform Bounded Distortion Lemma we assumed
that the
distortion of the composition depends on the sum of distortions of
individual
maps $h_{i}$. While that may often be a good estimate, it leaves the
case 
when distortions of maps cancel without satisfactory solution because
it does
not offer any way in which we could account for cancellations. For
example, high iterates of critical maps of the circle can be regarded as
compositions of the form considered by us with maps $h_{i}$ being
nearly 
linear-
fractional with nonlinearity totaling to close to zero. The basic
question is
whether in that situation we can ignore the maps $h_{i}$ completely
and 
approximate the whole composition ${\cal P}(f)$ by composition of maps
${\cal P}(\sigma_{i})$ only. The Cancellation Lemma formulated below
is a 
good tool to be used in such situations. The Lemma is a nice
illustration 
of the power of the Poincar\'{e}\  model approach. To the best of the
author's 
knowledge, no other technique has yielded a similar result. 

\paragraph{Cancellation Lemma.}
Let us consider a standard composition defined on an interval $(a,b)$
\[f=f_{m}=\sigma_{m}\circ H_{m}\circ\cdots\circ\sigma_{1}\circ H_{1}\;
.\]

We further assume that each map $H_{i}$ can be written as a
composition
\[H_{i}=h_{i}\circ g_{i}\]
where $h_{i}$ is a linear-fractional map.

We denote
\[\tilde{D}:=\sum_{j=1}^{k}\|{\cal P}(g_{j})-id\|_{C^{0}}\] and
\[\Delta:=\max\{|\sum_{i=1}^{j}({\cal P}(h_{j})-id)|: 1\leq j\leq
k\}\; .\]
\[S_{m}:={\cal P}(\sigma_{m})\ldots{\cal P}(\sigma_{1})\]

Then,
\[\|{\cal P}(f)-S_{m}(x)\|_{C^{0}}\leq\tilde{D}+2\Delta\; .\]

\subparagraph{A comment.}
What the Cancellation Lemma tells us is that if we have a composition
where 
maps of not 
neccessarily positive Schwarzian are all almost linear-fractional and
their 
distortions
almost cancel, then we can replace these maps with affine maps and
still get 
a good aproximation of the composition, at least locally. The main
value of this lemma is that its
assumptions are verified for first return maps of circle
homeomorphisms
as we prove in  the next section.

\paragraph{Beginning of the proof.}
The proof of the Cancellation Lemma is not very easy and will occupy
the next 
section. Here, we just make the first step which is the elimination of
maps
$g_{i}$ from the problem.

To achieve this, we put together maps $h_{i}$ and $\sigma_{i}$ and get
$s_{i}:=\sigma_{i}\circ h_{i}$. Then $f$ is a standard composition of
the
form
\[f=s_{m}\circ g_{m}\circ\cdot\circ s_{1}\circ g_{1}\; .\]

To this standard composition we apply Lemma ~\ref{lem1} and what we
get is
that 
\[{\cal P}(f)=G \circ {\cal P}(s_{m})\circ\cdots\circ {\cal
P}(s_{1})\]
in which
\[\|G-id\|_{C^{0}}\leq \tilde{D}\; .\]
\subparagraph{The reduced Cancellation Lemma.}
As the preceding argument shows, the functions $g_{i}$ can be omitted
from the composition defining $f$ and if we then prove that
\[\|{\cal P}(f)-S_{m}(x)\|_{C^{0}}\leq2\Delta\]
this will immediately imply the Cancellation Lemma. So we make this
additional 
assumption and call the resulting auxiliary theorem the ``reduced '' 
Cancellation Lemma. 

\subsection{Proof of the reduced Lemma.}
\paragraph{Extending functions $h_{i}$.}

The standard composition $f$ is now written as
\[f=f_{m}=\sigma_{m}\circ h_{m}\circ\cdots\circ\sigma_{1}\circ h_{1}\]
with all maps $h_{i}$ linear-fractional.

We then extend functions $h_{i}$ to one-parameter families $h^{t}_{i}$
defined
by
\[{\cal P}(h_{i}^{t})(x)=x+t\cdot({\cal P}(h_{i})(x) -x )\;\; 0\leq
t\leq 1\;
 .\]

In particular, $h_{i}^{0}=id$ and $h^{1}_{i}=h_{i}$.
The maps like $f_{i}^{t}$  are then defined in the obvious way.
\subparagraph{More important notations.} To simplify our future
equations we 
define:
\begin{equation}\label{equ:3,1}
 \Omega_{i}^{t}(x) :=
{\cal P}(h^{t}_{i})\circ {\cal P}(f_{i-1}^{t})(x)-{\cal
P}(\sigma_{i-1})
\circ\ldots\circ {\cal P}(\sigma_{1})(x)\; 
\end{equation}
and
\begin{equation}\label{equ:3,2}
\delta_{i}:={\cal P}(h_{i})(x)-x
\end{equation}
which quantity is independent of the choice of $x$ provided $h_{i}$ is
linear-fractional.

\subsection{Estimates}
With these notations we are ready to prove our basic lemma.
\begin{lem}   \label{lem11}
\[\frac{d\Omega_{i}^{t}}{dt}(x)\leq 2\Delta\]
for $1\leq i\leq m$ and any $x$.
\end{lem}
\begin{proof}
We compute the derivative in question:
\[\frac{d}{dt}\Omega_{i|x}^{t}=\sum_{j=1}^{i}\prod_{l=j}^{i-1}
\frac{d{\cal P}(\sigma_{l})}{dx}_{|{\cal P}(h^{t}_{l}
\circ f^{t}_{l}(x))}\cdot\frac{d{\cal P}(h^{t}_{j})}{dt}({\cal P}
(f_{j-1}^{t})(x))\]
\begin{equation}\label{7,2}
=\sum_{j=1}^{i}(\prod_{l=j}^{i-1}\frac{d{\cal P}(\sigma_{l})}{dx}
_{|{\cal P}(h^{t}_{l}\circ f^{t}_{l}(x))}) \cdot \delta_{j}\; 
\end{equation}
Here, we used the fact that the maps ${\cal P}(h^{t}_{j})$ are all
isometries,
thus their derivatives can be skipped in the product.
For simplicity, we will denote 
\[\prod_{l=j}^{i-1}\frac{d{\cal P}(\sigma_{l})}{dx}
_{|{\cal P}(h^{t}_{l}\circ f^{t}_{l}(x))}:=a_{j}(t)\; \]
Since maps $\sigma_{i}$ were assumed to have non-negative Schwarzian,
their
Poincar\'{e}\  models are weak contractions; thus $a_{j+1}(t)\geq
a_{j}(t)$.

Formula ~\ref{7,2} can then be rewritten as
\[\frac{d\Omega_{i}^{t}}{dt}(x)=\sum_{j=1}^{i} a_{j}\delta_{j}\; .\]
Now we use the famous Abel's series transformation to bound this
quantity:
\[\sum_{j=1}^{i} a_{j}\delta_{j}=\sum_{j=1}^{i}
a_{j}(\sum_{k=0}^{j}\delta_{k}-\sum_{k=0}^{j-1} \delta_{k})\]
where we adopted the convention $\delta_{0}=0$.

This can be further rewritten as
\[\sum_{j=1}^{i-1}((a_{j}-a_{j+1})\sum_{k=0}^{j}\delta_{k}) +
\sum_{k=0}^{i}
a_{i}\delta_{k}\; .\] 

The absolute value of the last expression is easy to bound. The last
term 
does not exceed $\Delta$ and the first one can be bounded by
\[ \sum_{j=1}^{i-1}|a_{j}-a_{j+1}|\Delta \leq \Delta\]
as the numbers $a_{j}$ form a non-decreasing sequence and are bounded
by $1$.

This concludes the proof.
\end{proof}

\paragraph{The conclusion.}

Lemma ~\ref{lem11} allows us to bound $\sup\{\Omega_{m}^{t}(x): x\in
R\}$  
by $2\Delta$, but this is exactly the statement of the reduced
Cancellation
Lemma.

\section{Pure Singularity Property}
A famous theorem of M. Herman says that smooth diffeomophisms of the
circle
are smoothly conjugated to rigid rotations, provided some
diophantine-type conditions are verified. One way to look at the smooth conjugacy is that in a small
scale it becomes $C^{1}$-close to linear. There is a dynamically
defined
event which takes place in a small scale. This is the first return map
to
a small interval. Hence, the smooth conjugacy tells us that first
return maps
for any diffeomorphism will tend to the first return maps for the
rigid rotation, that is, to linear maps. 

The first return map  is an ``induced map", which means that piecewise
it is a high iterate of the initial map. It is also known  that all intermediate images 
of the pieces of
its domain are disjoint and cover the circle completely. On each
piece, the
distortion is the total of distortions acquired by consecutive
iterates on corresponding intermediate images. Precisely, we can consider nonlinearity defined
in the introduction and then it turns out
that the nonlinearity of the first return map on each piece 
of its domain is equal to the sum of nonlinearities in all
intermediate
images transported by iterates of the map.

Since the linear map is characterized by ${\cal N}f=0$, in Herman's
theorem the 
distortions must
cancel. This not so surprising, perhaps, since the integral of the
nonlinearity
over the whole circle is $0$. But it is a remarkable fact  and in
certain approaches the central issue in the proof of Herman's theorem.

Naturally, we would like to know to what extent this fact is true for
critical
circle maps. For simplicity, we will consider maps with only one
critical point
such that coordinates can be changed $C^{3}$ smoothly to make it
locally
the map $x\rightarrow x^{\beta}+\epsilon$. It is widely conjectured,
and in few
very restricted cases has been argued with a computer's assistance,
that if we
choose
the domains of the first return map suitably , the sequence of the
first return maps will also approach a unique limit. But this limit map  is everything but 
 linear: its distortion does not not vanish and neither does its
Schwarzian 
derivative.  

Nevertheless, we prove that, in a certain sense, only distortions
acquired in an immediate proximity of the singularity count.
Distortions
due to the part remote from the singularity will tend to cancel, just
like
in the diffeomorphisms' case.

The proof that distortions cancel that we give is somewhat similar to
the
proof of Herman's Theorem given in \cite{Skh}.
\subsection{Assumptions and the statement of results}
\paragraph{The class of maps we are working with.} 
We will consider orientation-preserving $C^{3}$-smooth circle
homeorphisms
with one critical point of the polynomial type, at this moment of any
rotation number. 

As a consequence, (see \cite{mysa}), the circle is covered by two
overlapping 
open arcs. There is a ``remote" arc on which we assume
that the first derivative is bounded away from $0$. On the other
``close" arc 
the map has non-positive Schwarzian derivative.

We reserve the notation $f$ for maps in this class.

These assumptions are a little bit stronger than necessary for our
estimates 
to work, but we prefer not to obscure the idea by technicalities at
this 
point. 
We will discuss weakening of our requirements in the course of the
paper.

\paragraph{Some terminology.}
A {\bf symmetric neighborhood} is a neighborhood of the critical point
which
is contained in the close arc and the derivative of the function is
the same
in both endpoints. It follows from our assumptions that a symmetric 
neighborhood is also almost symmetric in the ordinary sense - the
critical 
point is in a bounded Poincar\'{e}\  distance from the the midpoint.

A {\bf chain of intervals} is a sequence of intervals such that each
is mapped
onto the next by the map. We will be particularly interested in chains
of
disjoint intervals. Obviously, there always is a map associated with a
chain,
namely the composition leading from the first interval to the last
one.

The continued fraction approximants of the rotation number will be
denoted with
\[\frac{p_{n}}{q_{n}}\; .\]
The denominators $q_{n}$ are important from the dynamical point of
view,
since they determine the times of closest returns by the orbit of a
point to
the point itself. The numbers satisfy the relations:
\[ q_{-1} = 0\; , q_{0} = 1\; , q_{n+1} = a_{n} q_{n} + q_{n-1}\]
where the coefficients $a_{n}$ are defined by the continued fraction
expansion
of the rotation number. An elementary discussion of the topological
dynamics
of diffeomorphisms with an irrational rotation number can be found in 
\cite{Skh}.  

An interval $J$ is said to be of the {\bf $j$-th order of fineness} if
\[j=\max\{i:\forall x\in J\;\, f^{q_{i}}(x)\notin J\}+1\; .\]

A {\bf uniform constant} is a function on our class of maps which
continuously
depends only on the quasisymmetric norm of the map, the logarithm of
the size 
of the close arc, the lower bound of the derivative on the remote arc,
and the 
$C^{3}$ norm. 

In view of this definition, ``absolute'' constants like $e^{\pi}$ are
uniform.
Perhaps a more meaningful example is the statement:

\begin{em} Let $f$ be a smooth circle homeomorphism with an irrational
diophantine rotation number. For each natural $n$, the derivative of
$f^{n}$ 
is bounded by a uniform constant.
\end{em}

Without the word ``uniform'' the statement would be  obviously true.
As it is
now, the main problem is whether the bound depends on $n$. If $f$ is a
diffeomorphism, it does not, and the sentence remains
true.\footnote{Which 
follows from Herman's theorem.} If $f$ is a critical map, the
statement is
false: the derivatives must become very large as $n$
grows.\footnote{Otherwise,
the map would be smoothly conjugated to the rotation, see
\cite{Herman}.}
 
\subparagraph{A notational convention.}
There will be so many uniform constants in use in the future
discussion that
we feel a need for a special notational convention to handle them
effectively.
Notations like  $K_{\cdot}$ will be used exclusively for uniform
constants.
The subscript will identify the particular constant.  All uniform
constants
will be introduced in lemmas, propositions or facts. The rule is that
in the
statement in which a constant is first defined, as well as in its
proof, the
constant will be identified by a single numerical subscript. The same
subscript
may denote different constants in different lemmas. 

However, when we use the constant later, its single subscript will be 
followed by an indication of where it was introduced. For example, the
constant
$K_{1}$ introduced in Lemma 10.15 will be called $K_{1}$ in the proof
of Lemma 
10.15, but later will be referred to as $K_{1,L.10.15}$.

\subparagraph{Approximate maps.}
To describe the procedure of approximating maps we need two objects :
a neighborhood\footnote{usually symmetric} of the critical point and a
chain of
intervals. To approximate the composition associated with this chain
we do the 
following:
\begin{itemize}
\item
If an interval is contained in the neighborhood, we leave the map
defined on it
to be $f$.
\item
Otherwise, instead of $f$ we use an affine map with the same image.
\end{itemize}

\paragraph{Main Theorem.}
\begin{em}
Let us suppose we have a chain of intervals 
\[(a_{0},b_{0}), (a_{1},b_{1}),\ldots,(a_{m},b_{m})\; ,\]
none of which contains the critical point, of the $\kappa$-th order of
fineness
and a symmetric neighborhood $U$ with the fineness of order $\lambda$.
This 
also
assumes that the length of the continued fraction expansion of the
rotation
number is at least $\kappa$.
We then approximate $f^{m}$ on $(a_{0},b_{0})$ to get some $\phi$. The
result
is that:

If $\kappa>\lambda$ then:
\[\|{\cal P}(f^{m})-{\cal P}(\phi)\|_{C^{0}}\leq
K_{1}K_{2}^{\sqrt{\kappa-\lambda}}\]
with $K_{1}$ and $K_{2}$ being uniform constants depending only on
global distortion properties of $f$ and $K_{2}<1$.
\end{em}

\subparagraph{A comment.}
The Main Theorem proves what we want to call informally ``the pure
singularity
property". If we have enough smoothness we can change coordinates 
so the resulting map is in our class and moreover its critical point
is locally 
in the form $x\rightarrow x^{\beta}+f(0)$. The Main Theorem then
asserts that 
asymptotically only what happens in this small neighborhood matters
and that is why
the expression ``pure singularity property" seems appropriate to the
author.

\subsection{General strategy of the proof}
The proof will largely use the concept of Poincar\'{e}\  model of the
interval
which is explained in earlier sections of this paper. We will also use
the same notations in this section. The main tool of our proof is going to be the cancellation Lemma.

\subparagraph{An outline of the argument.}
There is an obvious way we can use the Cancellation Lemma to prove our
Main 
Theorem.
Maps $\sigma_{i}$ will be the iterates of $f^{-1}$ which are left
unchanged 
when we replace $f^{-m}$ by its approximate 
map. Then, the claim of the Cancellation Lemma is exactly what we
assert in the theorem. Hence, our effort will be aimed towards verifying the hypotheses
of the Cancellation Lemma and then the main technical problem will be
to bound $\Delta$. We address these issues in the next section.

\section{Proof of the Main Theorem.}
We assume that the assumptions of the main theorem hold; in particular
that we
are given a chain of disjoint intervals. 
\subsection{Bounded distortion of critical circle maps}

\paragraph{Dynamical partitions.}
The forward orbit of the critical point $0$ defines a sequence of
partitions 
of the
circle, called dynamical partitions. For any $k$ less than the length
of the
continued fraction representation of the rotation number, the
$0,\ldots,a_{k}
q_{k}$
images of the interval $(0,f^{q_{k-1}}(0))$ are called {\bf lengthy
intervals}.
The $0,\ldots,q_{k-1}-1$ images of $(f^{q_{k}}(0),0)$ are called {\bf
short 
intervals}. Together, lengthy and short intervals form a partition of
the circle
and this is exactly what we are going to call the {\bf dynamical
partition of
the k-th order} and will be denoted $D_{k}$. 
\subparagraph{Consecutive dynamical partitions as refinements.}
We will now examine a very simple correspondence between dynamical
partitions
of order $k$ and $k+1$. The latter clearly is a refinement of the
former. More
precisely, all short intervals of the partition of order $k$ will
become lengthy
intervals of the next partition, while the lengthy intervals of the
coarser
partition will be subdivided. Each lengthy element of the partition of
order $k$
will be split into a number of lengthy intervals as well as one short
interval
which all belong to the partition of order $k+1$.

\subparagraph{Bounded geometry of dynamical partitions.}
The properties which are commonly referred to as ``bounded geometry''
are
summarized by the following statement:
\begin{fact}\label{fa:2,1}
If $f$ satisfies our regularity conditions, then:
\begin{itemize}
\item
The ratio of lengths of two adjacent elements of any dynamical
partition
is bounded by a uniform constant $K_{1}$.
\item
For any element of any dynamical partition, the ratios of its length
to the lengths of extreme intervals of next partition subdividing it 
are bounded by a constant $K_{2}$.
\end{itemize}
\end{fact}

Unfortunately, this Fact belongs to the ``folk
wisdom" and there is no clear reference to the proof. In one
particular
case when orbits of the critical point are periodic, Fact \ref{fa:2,1}
was verified in the work \cite{mysa}. Then, M. Herman showed how to 
carry this sort of estimates over to the more general situation (see
\cite
{michel}.)

\begin{lem}\label{lem:8,1}
There is a uniform constant $K_{1}$ such that the elements
of $D_{\lambda+K_{1}}$  adjacent to $0$ are contained in
$U$.\footnote{We 
remind the reader that $U$ is fixed and defined in the statement of
the 
Main Theorem.}
\end{lem}
\begin{proof}
This follows immediately from Fact ~\ref{fa:2,1} if we take into an
account
that $U$ is symmetric.
\end{proof}

\paragraph{Coarseness of dynamical partitions.}

We will need the fact that $f^{-1}$ on elements of $D_{j}$ with $j$
much 
larger than $\lambda$ is almost linear-fractional with almost constant
nonlinearity. To prove this, it is crucial to know that the partition
$D_{j}$ is fine enough. To establish this fact is the purpose of the
following
three lemmas.

\begin{defi}\label{defi:10,1}
Given a set $V$, the {\bf coarseness} of the partition $D_{j}$ {\bf
outside
$V$}, denoted $c_{j}(V)$ is defined by:
\[c_{j}(V)=\sum_{I\in D_{j}\: :\: I\cup(S^{1}- V)\neq \emptyset} 
(\frac{|I|}{dist\: (I,0)})^{2}\; .\]
\end{defi}

\begin{lem}\label{lem:10,1}
Let us fix $j$ and let $V$ be the interior of the union of two
elements of 
$D_{j}$ adjacent to $0$. Then, 
\[c_{j+1}(V) < K_{1} c_{j}(V)\]
where $K_{1}$ is a uniform constant less than $1$.
\end{lem}
\begin{proof}
This is an immediate corollary to Fact ~\ref{fa:2,1} and the
definition of
coarseness.
\end{proof}

\begin{lem}\label{lem:10,2}
Let $V$ and $j$ be related as in the statement of Lemma
~\ref{lem:10,1}.
Then $c_{j}(V)$ is uniformly bounded by some $K_{1}$.
\end{lem}
\begin{proof}
Let us define $V'$ to be the union of the elements of $D_{j+1}$
adjacent
to $0$. The first observation is that $c_{j+1}((S^{1}\setminus V)\cup
V')$
is uniformly bounded as a consequence of Fact ~\ref{fa:2,1}, second
part.
This in conjunction with Lemma ~\ref{lem:10,1} implies that
$c_{j+1}(V')$
can be bounded recursively by $c_{j}(V)$ times a constant less than
$1$
increased by a bounded amount. This recursive bound implies a bound
uniform in $j$.
\end{proof}

\begin{lem}\label{lem:10,3}
For $j>\lambda+K_{1,L.\ref{lem:8,1}}$, 
\[c_{j}(U)\leq K_{1}K_{2}^{j-\lambda}\]
with $K_{2}<1$.
\end{lem}
\begin{proof}
This follows immediately from Lemmas \ref{lem:8,1}, \ref{lem:10,1} and
\ref{lem:10,2}.
\end{proof}

\subsection{Preparations to use the Cancellation Lemma.}
\paragraph{How to represent $f^{-m}$ in the Cancellation Lemma ?}

We choose maps $\sigma_{i}$ to be the iterates of $f^{-1}$ on
intervals 
contained in $f(U)$ and others will be $H_{i}$. To obtain a
composition
in the form postulated by the hypotheses of the Cancellation Lemma we
may have
also to insert maps $\sigma_{i}$ equal to identities between
consecutive
maps $H_{i}$.

\subparagraph{Further choices.}  
The problem which still remains is a judicious choice of maps $h_{i}$.
We will simply give a prescription:
\begin{em}
We look at the nonlinearity of $H_{i}$ and $h_{i}$ will be the 
homography
that maps the domain of $H_{i}$ onto its image, and satisfies
\[{\cal P}(h_{i})(x)=x-1/2\cdot\int_{Dm(H_{i})}{\cal N}(H_{i})\; .\]
\end{em}
\paragraph{The assumptions of the Cancellation Lemma are then
satisfied.}
Next thing we need is to estimate constants which appear in the
Cancellation
Lemma.
It is relatively easy to deal with $\tilde{D}$ and we are going to
consider 
it first.
 
We introduce a map $g'_{i}$ as the map with the same image and
preimage as 
$H_{i}$, the same nonlinearity integral and constant nonlinearity.
First, we will estimate the Poincar\'{e}\  discrepancy between $H_{i}$
and $g'_{i}$ .

\begin{lem}\label{lem:2,1}
The difference
\[|{\cal P}(H_{i})-{\cal P}(g'_{i})|\]
is uniformly bounded by 
\[K_{1}(\frac{|(a_{i},b_{i})|}{dist\: ((a_{i},b_{i}),0)})^{2}\]
\end{lem}
\begin{proof}

First we precompose the maps with affine functions so as to have them
defined
on the unit interval. We observe that the nonlinearity coefficients
of the rescaled map are of the order of
\[ \frac{|(a_{i},b_{i})|}{dist\: ((a_{i},b_{i}),0)} \; .\]
With a slight abuse of notation we will still use the same symbols to
denote
the rescaled maps.

Let $\cdot^{*}$ denote the transport of forms by functions, and we
compute  
\[{\cal N}(H^{-1}_{i}\circ g'_{i})=(g'_{i})^{*}({\cal N}(H^{-1}_{i})+
{\cal N}(g'_{i})=\]
\[=(g'_{i})^{*}({\cal N}(H^{-1}_{i})+
{\cal N}(g'_{i})- (g'_{i})^{*}({\cal
N}((g')^{-1}_{i})+(g'_{i})^{*}({\cal N}((g')^{-1}_{i})=\]
\[=(g'_{i})^{*}({\cal N}(H^{-1}_{i}))-(g'_{i})^{*}({\cal
N}((g')^{-1}_{i}))\; .\]

Hence, it is enough to estimate the coefficient of 
\[{\cal N}(H^{-1}_{i})-{\cal N}(g')^{-1}_{i}=(H_{i}^{-1})^{*}{\cal
N}(H_{i})-((g')_{i}^{-1})^{*}{\cal N}(g'_{i})\; .\]

Finally introducing the derivatives explicitly we get:
\[(H_{i}^{-1})^{*}{\cal N}(H_{i})-(g_{i}^{,-1})^{*}{\cal N}(g'_{i})\]
\[={\cal N}(H_{i})(\frac{dH_{i}^{-1}}{dx}-\frac{d(g')_{i}^{-1}}{dx})+
\frac{d(g')_{i}^{-1}}{dx}({\cal N}(H_{i})-{\cal N}(g'_{i}))\; .\]

But it is evident that the estimate we want follows. Since the
coefficients of 
 nonlinearities are bounded as noted at the beginning of the proof,
the derivatives are 
\[1+O(\frac{|(a_{i},b_{i})|}{dist\: ((a_{i},b_{i}),0)})\; .\]
  
The lemma follows.
\end{proof}

\begin{lem}\label{lem:2,2}
The difference
\[|{\cal P}(h_{i})-{\cal P}(g'_{i})|\]
is uniformly bounded by 
\[K_{1}(\frac{|(a_{i},b_{i})|}{dist\: ((a_{i},b_{i}),0)})^{2}\; .\]
\end{lem}

\begin{proof}
By solving a corresponding differential equation, we find out that
a function with constant nonlinearity from the interval $[0,1]$ to
itself
is given by
\[x\rightarrow \frac{\exp(nx)-1}{\exp n -1}\]
where $n$ is the nonlinearity coefficient.

Next, we are going to find the displacement in the Poincar\'{e}\
model for
the image of any point $x\in (0,1)$. We will discard terms quadratic
or
higher in $n$.

If $x'$ means the image of $x$ we get
\[x'=x+\frac{1}{2}nx^{2}-\frac{1}{2}nx+O(n^{2})=x(1+\frac{1}{2}nx-\frac{1}{2}n
)+O(n^{2}) \]
and
\[1-x'=1-x-\frac{1}{2}nx^{2}+\frac{1}{2}nx+O(n^{2})=(1-x)(1+\frac{1}{2}nx)+
O(n^{2})\; .\]
The Poincar\'{e}\  displacement is 
\[-\log \frac{x(1-x')}{(1-x)x'}=-\log(
\frac{1+\frac{1}{2}nx}{1+\frac{1}{2}nx-
\frac{1}{2}n}+O(n^{2}))\]
\[=\log(\frac{1}{2}nx-\frac{1}{2}n-\frac{1}{2}nx+O(n^{2}))=-\frac{1}{2}n
+O(n^{2})\; .\]
Since $n$ is of the order of 
\[\frac{|(a_{i},b_{i})|}{dist\: ((a_{i},b_{i}),0)}\]
this concludes the proof.
\end{proof}

Our efforts are crowned by the following proposition:
\begin{prop}\label{prop:2,1}
If a chain of disjoint intervals satisfies the assumptions of Main
Theorem,
then our choice of maps $g_{i}$ and $h_{i}$ gives the bound for
$\tilde{D}$ by
\[\tilde{D}\leq K_{1}K_{2}^{\kappa-\lambda}\]
with $K_{2}<1$.
\end{prop}
\begin{proof}
Follows immediately from Lemmas ~\ref{lem:2,1}, \ref{lem:2,2} and 
\ref{lem:10,3}.
\end{proof}

\subparagraph{A technical comment.}
We could have done our estimates separately on the close arc and the
remote 
arc. While the author does not see any better method on the close arc,
on the remote arc it would be enough to assume that the second
derivative
is only H\"{o}lder continuous.

\paragraph{How to bound $\Delta$ ?}

Our main remaining problem is an estimate for $\Delta$. We have to be
able to 
see that this actually tends to $0$. By the definition of $g_{i}$, 
\begin{equation}\label{equ:5,1}
\sum_{i=0}^{j}({\cal P}(g_{i})(x)-x)=1/2\cdot \int_{C_{j}}{\cal N}f
\end{equation}
where $C_{j}$ is a union of these intervals of the chain from $0$ to
some $i_{j}$ that are not contained in $U$.

In order to prove Main Theorem we need to show that this integral is
uniformly 
exponentially small in $\sqrt{\kappa-\lambda}$. We call this the main
estimate.

\subsection{Main estimate}
\paragraph{Beginning of the proof of the main estimate.}
We consider the set $\tilde{U}$ which is equal to the complement of
$U$
together with intervals of the chain that are partly contained in $U$.
We regard $\tilde{U}$ with the suitably normalized Lebesgue measure as
a probability space. Then, the characteristic function $\chi$ of the
chain
can be viewed as a random variable, the dynamical partitions are
$\sigma$-
algebras of events, and individual elements of those partitions are
events. 

We will consequently use the language of conditional expectations,
denoted
with $E(\chi|\cdot)$ where the dot can be either a partition or a set
(event.)
The intuitive interpretation of $E(\chi|{\cal D}_{j})$, for example,
is the 
function whose value on each element of the partition is equal to the
relative
measure of the chain on this element.

The main estimate will follow from this proposition:

\begin{prop}\label{prop:2,2}
\[\int_{\tilde{U}}|E(\chi|{\cal D}_{j})(x)-E(\chi)|<K_{1} \cdot K_{2}^
{\sqrt{\kappa-j}} \]
for any $x$ in the element of $D_{j}$ completely disjoint with $U$
where $\lambda<j<\kappa$ and $K_{2}<1$, both $K_{1}$ and $K_{2}$ being
uniform 
constants
\footnote{We also remind the reader that $\kappa$ is the fineness 
of the chain.}.
\end{prop}

The proof of this proposition is lengthy and will occupy several next
pages.
 This is the main technical part of the proof of the main estimate,
and 
therefore the Main Theorem as well. 
\paragraph{The strategy of the proof.}
Let us fix a number $j$ as indicated in Proposition ~\ref{prop:2,2}.
In most arguments we will have to distinguish between two cases. 

If the chain contains only few intervals, then both $E(\chi)$ and 
$E(\chi|{\cal D}_{j})$ for $j$ larger than $\lambda$ are small. In
particular,
we will assume without loss of generality that $E(\chi|{\cal D}_{j}) >
0$.
Otherwise, there is an element of ${\cal D}_{j}$ with no intervals
from
the chain in it. Then, any other element of that partition may
intersect with
at most one interval. Thus, $E(\chi|{\cal D}_{j})$ would be
exponentially 
small  with respect to $\kappa-j$ by bounded geometry (see Fact
~\ref{fa:2,1}.)

The second, more important and 
interesting situation is when the number of intervals is large, and
individual 
expectations of $\chi$ are much greater than our estimate. Then we
will 
use ``averaging'' techniques based on the fact that the nonlinearity
integral
over the whole $\tilde{U}$ is close to $0$. Unfortunately, these
approaches
are mutually exclusive: we cannot use averaging if we only have few
intervals
in the chain. Thus, we will have to keep track of both possibilities 
throughout ther proof of Proposition ~\ref{prop:2,2}. 

We will look at the quantity
\[v_{j}=\frac{\max\{E(\chi|{\cal D}_{j})(x):x\notin
U\}}{\min\{E(\chi|{\cal D}_{j})(x):x\notin U\}}-1\; .\]
We will prove that if $j$ is more than $\lambda$ but sufficiently
smaller
than $\kappa$, this quantity will show definite growth as $j$ is
increased.
On the other hand, we will notice that there are certain bounds on its
growth
and this will let us assert that its initial value $v_{\lambda}$ must
be very 
small. Proposition ~\ref{prop:2,2} will follow in this way.

The conditional expectations $E(\chi|D_{\cdot})$ will be referred to
shortly
as {\bf densities}.

\paragraph{Technical preparations.}
\subparagraph{Approximate invariance of $\chi$.}
\begin{lem}\label{lem:5,3}
Let $J_{1}$ and $J_{2}$ be two similar\footnote{That is, both lengthy
or both
short.} elements of $D_{j}$. We further assume that neither of them is
contained in $U$. To fix notations, we assume that
$J_{1}=f^{l}(J_{2})$ where
$l$ may be negative, but $|l|<q_{j}$.

If 
\[\kappa-K_{1}|\log(E(\chi|D_{j})(J_{1}))|>j\; ,\]
then
\[|\frac{E(\chi|D_{j})(J_{1})}{E(\chi\circ f^{l}|D_{j})(J_{1})}-1|<
K_{2}\frac{\exp(j-\kappa)}{E(\chi|J_{1})}\; .\]
\end{lem}
\begin{proof}
Since $\chi$ is the characteristic function of a chain, $\chi$ and 
$\chi\circ f$ may differ on at most two intervals of the size order
$\kappa$.
By bounded geometry, the length of an interval of fineness $\kappa$
 is exponentially small with the exponent $\kappa-j$ compared to the
element of
$D_{j}$ that contains it. So first we pick $K_{1}$ to ensure that
$J_{1}$ 
contains at least two intervals of the chain. 

If this is true, the left-hand side of the second condition is roughly
the ratio
of the measure of one interval from the chain to the total measure of
the 
portion of the chain in the
the $j$-th partition containing it. This means that
\[|\frac{E(\chi|D_{j})(J_{1})}{E(\chi\circ f^{l}|D_{j})(J_{1})}-1|<
K_{2}\frac{\exp(j-\kappa)\cdot |J_{1}|}{\int_{J_{1}}\chi}\] which
immediately 
yields the claim of the lemma.

\end{proof} 
\subparagraph{Bounds on $v_{i}$.}
\begin{lem}\label{lem:5,2}
If $\kappa-K_{2}|\log E(\chi)|>j>\lambda+K_{1}$ then $v_{j}<K_{3}$.
\end{lem}
\begin{proof}
There must be an interval $J\in D_{j}$ which is not adjacent to $0$
with density
comparable to $\int\chi$. Then by Lemma \ref{lem:5,3} we see it is
possible to
find $K_{2}$ so that $\chi$ is sufficiently close to $\chi\circ f$ and
since
by choosing $K_{1}$ we can make sure that
the jacobian $f'$ has bounded variation when we go to an interval
similar to $J$, we immediately obtain that the densities on intervals similar to $J$ are
all comparable to $\int\chi$. 

Let us  assume that $J$ is a lengthy interval.
What we still need to show is the density on short intervals is also
at least
comparable. It cannot be much larger, since the images of short
intervals will
cover a definite portion of $J$ in the next subdivision.  

\begin{itemize}
\item
The subdivision of $J$ contains a lengthy interval of $D_{j+1}$ with
density
at least comparable to the density on $J$. Then, since this lengthy
interval is
an image of any short interval from $D_{j}$, the estimate follows.
\item
All lengthy intervals subdividing $J$ have very small density compared
with
the density on $J$. Then, the remaining short interval of $D_{j+1}$
must 
attain very high density. But the images of this interval will fill up
a definite portion of short intervals of $D_{j}$ and again the
estimate 
follows.
\end{itemize}  

The case when $J$ is a short interval can be solved using a similar,
but 
easier, reasoning.
\end{proof}

\subparagraph{Another technical lemma.}
\begin{lem}\label{lem:2,5}
Suppose that we have two functions $g$ and $h$ on the unit interval
and both
 positive and measurable
with respect to the same $\sigma$-algebra. Let us also assume that
$\int h=1$
and $\inf\{h(x) : x\in [0,1]\}=C_{1}>0$. Suppose
further that $\int g = 1$ and $\int gh = 1+\epsilon\, ,\epsilon>0$.
Then,
\[\frac{\max\{g(x):x\in I\}}{\min\{g(x):x\in
I\}}>1+C_{2}(C_{1})\epsilon\; \]
where $C_{2}$ is a continuous function of $C_{1}$ only.
\end{lem}
\begin{proof}
We introduce
\[g':=g-1-\epsilon/2\; .\]
We also define $I_{+}=\{x\in I:g'\geq 0\}$ and $I_{-}=I\setminus
I_{+}$.
We have
\[\int_{I_{+}} g'h + \int_{I_{-}} g`h = \epsilon/2\; .\]
But 
\[\int_{I_{-}} g' <- \epsilon/2\]
and thus 
\[\int_{I_{-}} hg'< - C_{1}\epsilon/2\]
Hence
\[\max\{g'(x):x\in I\}>\int_{I_{+}} g'h > \epsilon/2(1+C_{1})\; .\]
\end{proof}

\paragraph{The minimum and maximum density in similar intervals.}
\begin{lem}\label{lem:2,4}
Let us choose an integer $\hat{j}$ which satisfies
\[\lambda+K_{1,L.\ref{lem:5,2}}<\hat{j}<\kappa\; .\]
Suppose that the smallest density is attained on an 
element $Y$ of ${\cal D}_{\hat{j}}$. and the largest density is on
$X\in {\cal D}_{\hat{j}}$.

Then, there are three possibilities:
\begin{enumerate}
\item
\[K_{2}\frac{\exp(\hat{j}-\kappa)}{E(\chi|D_{j})(X)}>K_{3}v_{\hat{j}}\;
.\]
\item
\[E(\chi)<K_{6}K_{7}^{\kappa-\hat{j}}\]
where $K_{7}<1$
\item
\[v_{\hat{j}+[K_{4}|\log(v_{\hat{j}})|]}-v_{\hat{j}}>K_{5}v_{\hat{j}}\]
\end{enumerate}

All constants mentioned are positive.\footnote{And the square brackets
in the third condition mean the ``ceil'' function:
\[[x] := \inf\{n\in {\bf Z} : n\geq x\}\; .\]}
\end{lem}
\begin{proof}
We first explain the role of $K_{4}$.
By the bounded geometry and Koebe principle, the iterate of $f$ which
maps $X$
to $Y$,denoted $\phi$, has bounded nonlinearity on $X$. By the bounded
geometry again, the
maximum length of elements of ${\cal D}_{\hat{j}+l}$ on  $X$ tends to
$0$ 
exponentially fast in $l$.\footnote{The quickest way to see this is by
the
Uniform Bounded Distortion Lemma. But if one works a little harder and
proves 
that the Schwarzian has a definite sign for all but very low iterates,
the 
Koebe Priciple will be enough.} Therefore, $K_{4}$ can be chosen so
that for 
$l>K_{4}|\log(v_{\hat{j}})|$ the 
supremum of the quantity
$\phi'(x)-\phi'(y)$ on any element of ${\cal D}_{\hat{j}+l}$ is less 
than 
$K_{8}v_{\hat{j}}$
where $K_{8}$ can be made arbitraly small by choosing $K_{4}$
sufficiently
large.

Next, we will analyze the meaning of the first alternative.  When we
choose 
\[K_{2}=K_{2,L.\ref{lem:5,3}}\; .\]
In accordance with Lemma \ref{lem:5,3} the left-hand side of the
inequality
which is the first alternative bounds
\[\max(\frac{E(\chi|D_{\hat{j}})(X)}{E(\chi\circ\phi|D_{\hat{j}})(X)},\frac
{E(\chi\circ\phi|D_{\hat{j}}(X)}{E(\chi|D_{\hat{j}})(X)})-1\; .\]

By choosing $K_{3}$ appropriately, we can ensure that if the first
alternative
does not occur, the above quantity is as small compared with
$v_{\hat{j}}$ as
we may wish. 

Thus, we get
\begin{equation}\label{equ:6,1}
\int_{Y}\chi\cdot\phi'\geq (1-K_{3}v_{\hat{j}})\int_{X}\chi\; .
\end{equation}

So we assume that the first possibility in Lemma \ref{lem:2,4} does
not
occur which technically means Equation \ref{equ:6,1}. The appropriate
value
of $K_{3}$ will be chosen later. We also pick
$l=[K_{4}v_{\hat{j}}]+1$.

Then we have to think  what the second alternative means. According to
Lemma ~\ref{lem:5,2} either $v_{\hat{j}}<K_{3,L.\ref{lem:5,2}}$ or
$\int\chi$ is exponentially small
in $\kappa-\hat{j}$. Therefore, we can pick constants
$K_{6}$ and $K_{7}$ in such a way that if the second alternative does
not 
occur, then $v_{\hat{j}}$ is bounded by $K_{3,L.\ref{lem:5,2}}$.

So, we now assume that neither the first nor the second alternative
occurs,
which technically means that formula ~\ref{equ:6,1} holds with a small
$K_{3}$ 
to be arbitrarily chosen, and that $v_{j}$ is uniformly bounded.

By our assumption and Lemma ~\ref{lem:5,3}
\begin{equation}\label{equ:2,5}
\int_{Y}\chi\cdot \phi'\geq
(1+v_{\hat{j}}-K_{3}v_{\hat{j}}-K_{3}v^{2}_{\hat{j}})\int_{Y}\chi\; .
\end{equation}

Next, we condition $\chi$ with respect to $D_{\hat{j}+l}$:
\begin{equation}\label{equ:6,2}
\int_{Y}\chi\phi'=\int_{Y}\chi(\phi'-E(\phi'|D_{\hat{j}+l}
))+\int_{Y}
E(\chi|D_{\hat{j}+l})\cdot E(\phi'|D_{\hat{j}+l})
\end{equation}

The first term in Equation ~\ref{equ:6,2} is bounded by
\begin{equation}\label{equ:6,3}
\int_{Y}\chi(\phi'-E(\phi'|D_{\hat{j}+l})<K_{8}v_{\hat{j}}\int_{Y}E(\chi|D_{\hat{j}+l})
\end{equation}

Here, the reader may recall that $K_{8}$ is an auxiliary constant
which can be
made as small as we want by adjusting $K_{4}$.

Putting Equations \ref{equ:2,5}, \ref{equ:6,2} and \ref{equ:6,3}
together we
get
\begin{equation}\label{equ:2,6}
\int_{Y}E(\chi|D_{\hat{j}+l})\cdot E(\phi'|D_{\hat{j}+l})\geq
(1+v_{\hat{j}}-K_{8}v_{\hat{j}}-
K_{3}(1+v_{\hat{j}})v_{\hat{j}})\int_{Y}E(\chi|D_{\hat{j}+l})
\end{equation}

We note here that $v_{\hat{j}}$ is bounded by $K_{3,L.\ref{lem:5,2}}$
and so
we can pick $K_{3}$ and $K_{4}$ which controls $K_{8}$ in such a way
that 
\[ K_{5}=(1-K_{8}-K_{3}(1+v_{\hat{j}}))\cdot
K_{3,L.\ref{lem:2,5}}-1>0\; .\]

Our final step is to conclude the third alternative from this
inequality by 
Lemma ~\ref{lem:2,5}.
\end{proof}

\paragraph{Maximum and minimum density in intervals of different
kind.}
This case will be solved by the following lemma:
\begin{lem}\label{lem:2,6}
Let us suppose that the minimum density is attained on a short
interval $X$ and
the maximum density on a lengthy interval $Y$ and
$v_{\hat{j}}=\epsilon$. Then, we 
look at the next dynamical partition, in which $X$ becomes a lengthy
interval and $Y$ gets subdivided into a number of lengthy
intervals  $Y_{i}$ and one short $X'$. At least one of the following
holds 
true:
\begin{itemize}
\item
For some $i$ 
\[E(\chi|Y_{i})>K_{1}\cdot E(\chi|Y)\] with $K_{1}<1$, but, on the
other hand,
subject to the condition $K_{1}\cdot K_{5,L.\ref{lem:2,4}} > 1$.
\item 
\[v_{\hat{j}+1}>1+K_{2}\epsilon\] where $K_{2}>1$ is a uniform
constant.
\end{itemize}
\end{lem}
\begin{proof}
This is a rather obvious fact. We will not give a detailed proof.
Instead, we 
note informally that if the possibility does not occur it means that
on all 
intervals
$Y_{i}$ the density is smaller by a definite amount than the density
on the
whole $Y$. Then, the density on $X'$ has to exceed the density on $Y$.
Moreover, it must be larger by a definite amount, since the relative
measure of
$X'$ with respect to $Y$ is not too big by the bounded geometry.
\end{proof}

\subparagraph{A remark.}
Of course, the claim of Lemma ~\ref{lem:2,6} also holds if the largest
density 
is attained on a short interval and the smallest on a lengthy one, and
the 
proof is the same.

\paragraph{Proof of Proposition  ~\ref{prop:2,2}.}
With all these technical facts we are ready to conclude the proof of
Proposition
~\ref{prop:2,2}. 

Let us discuss the meaning of Lemmas ~\ref{lem:2,4} and
~\ref{lem:2,6}. 
Consider a $j\leq \hat{j}\leq (\kappa+j)/2$.

Suppose the maximum ratio of densities is attained on intervals of 
${\cal D}_{\hat{j}}$ which are of different kind. 

Then we want to use Lemma
~\ref{lem:2,4}, and suppose first that the first alternative in Lemma
occurs.
If we multiply the corresponding inequality on both sides by
$E(\chi|X)$ and
integrate over the whole $\tilde{U}$, we get what Proposition
~\ref{prop:2,2} 
claims, and even more, since there is no radical in the exponent.
However, 
we get this claim for $\hat{j}$, not $j$. Now, as
$j$ is no larger than $\hat{j}$, the left-hand side in Proposition
~\ref{prop:2,2} will not increase if we replace $\hat{j}$ with $j$. On the other hand, as
we assumed that $\hat{j} \leq (j+\kappa)/2$, the exponential
right-hand side
will suffer only a bouned decrease if make such a replacement. Thus, 
Proposition ~\ref{prop:2,2} is proven in that case.  

We leave it to the reader to show that Proposition ~\ref{prop:2,2}
also holds
if the second alternative occurs. 
  
If the third alternative in Lemma ~\ref{lem:2,4} is the only one, then
we want to replace $\hat{j}$ with 
$\hat{j}+[K_{4,P.\ref{lem:2,4}}\log v_{\hat{j}}]$. According to the
claim of 
the Lemma, $v_{\hat{j}}$ is bound to increase by a definite factor.

Now, let us think about Lemma ~\ref
{lem:2,6}. It says that in the next dynamical partition
$v_{\hat{j}+1}$ could 
be substantially more, or it may be less, but then the extrema are
attained on
intervals of the same kind. Moreover, in that second case the
constants have
been set up so that the increase implied by the subsequent use of
Lemma
~\ref{lem:2,4} will offset the tiny decrease on the previous step.

Thus, if start with $\hat{j}=j$ and follow this reasoning, we see that
either 
Proposition ~\ref{prop:2,2} holds or $v_{\hat{j}}$ grows
exponentially at a rate inversely proportional to
$|\log(v_{\hat{j}})|$ as
we increase $\hat{j}$. This means that as $\hat{j}$ finally reaches
$(\kappa+j)/2$, the $\log(v_{\hat{j}})$ will have grown by roughly 
$\sqrt{(\kappa-j)/2}$. But, in view of Lemma ~\ref{lem:5,2},
$v_{\hat{j}}$
is either uniformly bounded, or Proposition ~\ref{prop:2,2} holds
anyway, since
$E(\chi)$ is very small. Therefore, the initial value of $\log
v_{\hat{j}}$ 
must be of the order of $\sqrt{\kappa-j}$.

This concludes the proof of Proposition ~\ref{prop:2,2}.

\paragraph{Main Theorem follows easily.}

We wish to bound
\begin{equation}\label{equ:6,4}
\int_{\tilde{U}} \chi\cdot n 
\end{equation}
where $n$ is the nonlinearity coefficient.

We compute
\begin{equation}\label{equ:6,5}
 \int_{\tilde{U}}\chi n  = \int_{\tilde{U}}\chi (n-E(n|D_{j})) +
\int_{\tilde{U}} \chi E(n|D_{j}) 
\end{equation}

Here, $j$ is chosen so it satisfies the assumptions of Proposition
~\ref{prop:2,2}. The first term in Equation \ref{equ:6,5} is easily estimated using
Lemma ~\ref{lem:10,3}. The lemma implies that $n-E(n|D_{j})$ is
exponentially
small with $\kappa-\lambda$ and certainly the same is true of the
integral.

Our further effort is then aimed at estimating the second term.
\begin{equation}\label{equ:6,6}
\int_{\tilde{U}}\chi E(n|D_{j})  = \int_{\tilde{U}} E(\chi|D_{j})
E(n|D_{j}) 
\end{equation}
\[= \int_{\tilde{U}}(E(\chi|D_{j})-E(\chi)) E(n|D_{j})  +
\int_{\tilde{U}} E(n|D_{j})\; .\]

Again, Equation \ref{equ:6,6} reduces the problem to bounding two
terms.
The first term is bounded based on Proposition \ref{prop:2,2} since
the 
integral 
\[|E(\chi) - E(\chi|D_{j})|\]
is small and the nonlinearity coefficient is bounded.

In the second term the constant can be moved from the integral and the
remaining integral is exponentially close to $0$ as $U$ was assumed to be symmetric 
and $\tilde{U}$ differs from $U$ at most by a length of the interval
from
the chain.

This concludes the proof of the main estimate. Then, we can bound
$\Delta$
in the Cancellation Lemma which immediately yields our Main Theorem.

\subsection{Final remarks.}
\paragraph{The pure singularity property and unimodal maps.}
The pure singularity property for critical circle maps has an analogue
in the study of unimodal maps with the dynamics of solenoidal type. Since in such a
situation  the measure of the corresponding Cantor attractor is $0$
(see
\cite{blyu}), it follows immediately (for example by Lemma
~\ref{lem1}, see
also \cite{miszczu}) that the joint distortion of the first return map
due to 
parts 
of the interval where the nonlinearity is bounded must tend to $0$. So
the
counterpart of the pure singularity property also holds for unimodal
maps
and is not a very difficult fact in that context. Thus, the pure
singularity
property will allow us to extend certain results proved for unimodal
maps 
to circle maps. For example, we can follow the line of argument by W.
Pa\l uba 
which was originally developed in the context of unimodal maps with
solenoidal 
dynamics.
\footnote{I learnt about the result by personal communication and it
is to
be part of the his Ph.D. thesis.}  

The result we get for circle maps is this:

\begin{em}
Whenever two circle maps, each having a critical point of the same
type, are 
Lipschitz-conjugate, the conjugacy is differentiable at the critical
point.
\end{em}

\paragraph{Asymptotic analyticity?}
For a large class of circle maps we can change the coordinates in such
a way 
that the map becomes a polynomial in the neighborhood of the critical
point. 
One is lead to expect that since asymptotically the first return map
tends to 
the 
composition of pieces of polynomial maps. If this convergence could be
proven 
to be uniform in the vicinity of the circle on the complex plane, that
might be
an important step, possibly allowing the use of some machinery
developed in 
~\cite{miszczu}.
Unfortunately, the pure singularity property itself does not seem to
allow that
conclusion, as it gives the estimate on the circle only. Hopefully,
the 
conjecture will be proven in a forthcoming paper.


\newpage
\begin{thebibliography}{9}
\bibitem{miszczu} D. Sullivan: {\em On the structure of infinitely
many
dynamical systems nested inside or outside a given one}, preprint
IHES/M/90/75
\bibitem{preston} Preston,C.: {\em Iterates of maps on an interval},
Lecture Notes 
in Mathematics, Vol. 999. Berlin,Heidelberg,New York: Springer (1983)
\bibitem{Rand} D. Rand: {\em Global phase-space universality, smooth 
conjugacies and renormalisation I. The $C^{1+\alpha}$ case.},
Nonlinearity
{\bf 1} (1988) 181-202
\bibitem{dMvS} W. de Melo: {\em Lectures on one-dimensional dynamics},
17. 
Col\`{o}quio Brasileiro de Matem\`{a}tica, IMPA (1989)
\bibitem{guj} J. Guckenheimer: {\em Limit Sets of S-Unimodal Maps
with Zero Entropy}, Commun. Math. Phys. {\bf 110}, 655-659(1987)
\bibitem{blyu} A. Blokh, M. Lyubich: {\em Measure and dimension of 
solenoidal attractors of one-dimensional dynamical systems}, to appear
in 
Comm. in Math. Phys.
\bibitem{mysa} G. Swiatek: {\em Rational rotation numbers for maps of
the
circle}, Commun. in Math. Phys., {\bf 119}, 109-128 (1988) 
\bibitem{michel} M. Herman: {\em Conjugaison quasi sym\'{e}trique des 
hom\'{e}omorphismes analitique de cercle \`{a} des rotations}, a
manuscript
\bibitem{Skh} Ya. Sinai, K. Khanin: {\em A new proof of M. Herman's
theorem},
Commun. Math. Phys. {\bf 112}, 89-101 (1987)
\bibitem{Herman} M. Herman: {\em Sur la conjugaison differentiable des
diffeomorphismes du cercle a des Rotations}, Publ. Math. IHES {\bf
49}, (1976).
\bibitem{tanver} J.J.P. Veerman, F. Tangerman: {\em Scalings in
circlemaps I},
SUNY at Stony Brook IMS preprint no.8(1990), to appear in Commun.
Math. Phys. 

\end{thebibliography}
\end{document}